\documentclass[11pt,reqno, a4paper]{amsart}

\usepackage{amsmath,amsthm}
\usepackage{amssymb}
\usepackage{amsfonts}
\usepackage{amscd}
\usepackage{enumerate}
\usepackage{color}
 at 9truept

\newtheorem{thm}{Theorem}[section]
\newtheorem{theorem}[thm]{Theorem}

\newtheorem{corollary}[thm]{Corollary}
\newtheorem{lemma}[thm]{Lemma}
\newtheorem{proposition}[thm]{Proposition}
\newtheorem{definition}[thm]{Definition}

\theoremstyle{definition}
\newtheorem{example}[thm]{Example}
\newtheorem{examples}[thm]{Examples}
\newtheorem{remark}[thm]{Remark}
\newtheorem{remarks}[thm]{Remarks}

\newtheorem{notation}[thm]{Notation}
\newtheorem{free text}[thm]{}

\numberwithin{equation}{section}

\newcommand{\N} {\mathbb{N}}

\newcommand{\rmref}[1]{{\rm (}\ref{#1}{\rm )}}

\newcommand{\Hom}{\operatorname{Hom}}

\newcommand{\Mod}{\operatorname{Mod}}

\newcommand{\Tor}{{\rm Tor}}
\newcommand{\Ext}{{\rm Ext}}

\newcommand{\ga}{\alpha}

\newcommand{{\Uop}}{{U^{\rm op}}}
\newcommand{{\Ae}}{{A^{\rm e}}}
\newcommand{{\Aop}}{{A^{\rm op}}}
\newcommand{{\Bop}}{{B^{\rm op}}}

\newcommand{\ahha}{{\scriptscriptstyle{A}}}
\newcommand{\Aopp}{{\scriptscriptstyle{\Aop}}}

\newcommand{{\op}}{{{\rm op}}}
\newcommand{{\coop}}{{{\rm coop}}}
\newcommand{{\sop}}{{*^{\rm op}}}

\newcommand{\lact}{{\,\raise1pt\hbox{$\scriptscriptstyle{\rhd}$}\,}}                  %
\newcommand{\ract}{{\,\raise1pt\hbox{$\scriptscriptstyle{\lhd}$}\,}}                  
\newcommand{\blact}{{\,\raise1pt\hbox{$\scriptscriptstyle{\blacktriangleright}$}\,}}  %
\newcommand{\bract}{{\,\raise1pt\hbox{$\scriptscriptstyle{\blacktriangleleft}$}\,}}   %

\newcommand{\due}[3]{{}_{{#2 }} {#1}_{{ #3}}\,}    

\begin{document}

\title{Poincaré Duality for Hopf algebroids}
\author{Sophie Chemla}

\maketitle

\begin{abstract}
We prove a twisted Poincaré duality for (full) Hopf algebroids with bijective antipode. As an application, we recover the Hochschild twisted Poincaré duality of Van Den Bergh
 (\cite{VDB}). We also get a  Poisson twisted Poincaré duality, which was already stated for oriented Poisson manifolds  in \cite{CLYZ}. 
\end{abstract}
\maketitle

\section{Introduction}

Left bialgebroids over a (possibly) non commutative basis $A$ generalize bialgebras.
If $U$ is a left bialgebroid, there is a natural $U$-module structure on $A$ and 
the category of left modules on a left bialgebroid $U$ is monoidal. Nevertheless, $A$ is not necessarily a right $U$-module.
 Left Hopf left bialgebroids (or $\times_A$ -Hopf algebras)
\cite{Sch}) generalize Hopf algebras.   
 In a left Hopf left bialgebroid $U$,  the existence of an antipode is not required but, for any element $u \in U$, there exists an element $u_+\otimes u_- $ corresponding to $u_{(1)}\otimes S(u_{(2)})$.
 The structure of full Hopf algebroids  (\cite{BoeSzl}) ensures the existence of an antipode. 
  If $L$ is a Lie Rinehart algebra (or Lie algebroid) \cite{Rin}, there exists a standard left bialgebroid structure on its enveloping algebra $V(L)$. This structure is left Hopf. Kowalzig showed (\cite{Kow}) that $V(L)$ is a Hopf algebroid if and only if there exists a right $V(L)$-module structure on $A$.  If $X$ is a ${\mathcal C}^\infty$ Poisson manifold and $A={\mathcal C}^\infty (X)$,  the $A$-module of  global differential one forms $\Omega^1(X)$ is endowed with a natural Lie Rinehart  structure, which is of much interest. Moreover,  Huebschmann (\cite{Hue})
exhibited a right $V(\Omega^1(X))$-module structure on $A$ (denoted $A_P$) that makes $V(\Omega^1(X))$ a full Hopf algebroid. 
He 	also  interpreted  the Lichnerowicz Poisson cohomology $H^i_{Pois}(X) $  as 
$\Ext^i _{V(\Omega^1(X))}(A,A)$ and the Poisson homology $\Tor^{Pois}_i (X)$ (\cite{Bry}, \cite{Kos}) of $X$ as $\Tor_{V(\Omega^1(X))}(A_P, A)$. 

A Poincaré duality theorem was proved in \cite{Che} for Lie Rinehart algebras and then extended to left Hopf left bialgebroids in \cite{KowKra}. It asserts, under some conditions,  that if $\Ext^i_U(A,U)=0$ for $i\neq d$,  there is an isomorphism for all left $U$ module $M$ and all $n \in \N$ 
$$\Ext_U^n(A,M)\simeq \Tor^U _{d-n}(M \otimes \Lambda, A)$$
where $\Lambda := \Ext^d_U(A,U)$ with the right $U$-module structure given by right multiplication.  
If $U=V(L)$ is the enveloping algebra of a finitely generated projective  Lie Rinehart algebra $L$, it is shown in \cite{Che} that 
$\Ext^n_{V(L)}(A,V(L))=0$ if $n\neq dim L$. Moreover   $\Ext^{dim \; L}_{V(L)}(A,V(L))\simeq  \Lambda_A^{dim\; L}(L^*)$.

We give a new formulation of Poincaré duality in the case where $U$ is left Hopf,   its coopposite $U_{coop}$ as well and  $A$ is endowed with a right $U$-module structure. (denoted $A_R$). \\

{\bf Theorem} \ref{Poincare duality}
{\it  Let $U$ be a left and right Hopf left bialgebroid over $A$. Assume that 
\begin{enumerate}
\item $\Ext^{i}_U(A,U)=\{0\}\; if\;  i\neq d$ and set $\Lambda =\Ext^{d}_U(A,U)$.  
\item The left $U$-module $A$ admits a finitely generated projective resolution of finite length
\item $A$ is endowed with a right $U$-module structure (denoted $A_R$)  for which it is a dualizing module 
\item Let ${\mathcal T}$ be the left $U$-module $Hom_A ({A_R }_\bract , \Lambda_\bract )$ (see proposition \ref{structures}). The $A$ module $_\lact {\mathcal T}$ and the 
$A^{op}$-module $ {\mathcal T}_\ract$ are projective
\end{enumerate}
Then, for all $i \in \N$ there is an isomorphism
$$\Ext^i_{U}(A,M)\simeq \Tor^{U}_{d-i} \left ( A_R,  {\mathcal T}_\ract {\otimes_A} _\lact M \right ).$$ }

Assume now that $H$ is a full Hopf algebroid. The antipode allows to transform any left (resp. right) H-module $M$ (resp. $N$) into a right (resp. left ) $H$-module 
denoted $M_S$ (resp. $_S N$). Thus from the left $H$-module structure on $A$, we can make a right $H$-module structure $A_S$. From the right $H$-module structure on $\Lambda$, we can make a left $H$-module structure denoted $_S \Lambda$.  The duality states as follows :
$$\Ext^i_{H}(A,M)\simeq \Tor^{H}_{d-i} \left ( A_S,  {_S\Lambda} \otimes_A M \right ).$$ 
In the special case of the (full) Hopf algebroid $A\otimes A^{op}$, we recover the Hochschild twisted Poincaré duality of \cite{VDB}. 
In the special case where $X$ is a Poisson manifold and $H=V(\Omega^1(X))$, the duality above can be rewritten in terms of Poisson cohomology and 
homology. 
Let $M$ be a left $H$-module. The coproduct on $H$ allows to endow $_S\Lambda \otimes_A M$ with a left $H$-module structure. Denote by $H^i _{Pois}(M)$ the Poisson cohomology with coefficients in $M$ and 
$H_i^{Pois} \left ( {_S\Lambda} \otimes_A M \right )$ the Poisson homology with coefficients in $_S\Lambda \otimes_A M$. There is an isomorphism 
$$H^i _{Pois}(M)\simeq  H_{d-i}^{Pois} \left ( _S\Lambda \otimes_A M \right ).$$
This formula was stated in  \cite{CLYZ} for oriented Poisson manifolds (see also \cite{LauRic} for polynomial algebras with quadratic Poisson structures, \cite{Zhu} for linear Poisson structures, \cite{LuoWangWu} for general polynomial Poisson algebras).

\vspace{1em}

\centerline{\bf Acknowledgement}

\vspace{1em}
I thank N. Kowalzig for comments on an earlier version of this text. 

 \vspace{1em}
\centerline{\bf Notations}
\vspace{1em}
Fix an (associative, unital, commutative) ground ring  $ k $.  Unadorned tensor products will always be meant over  $ k $.
All other algebras, modules etc.~will
have an underlying structure of a
$k$-module. Secondly, fix an associative and unital
$k$-algebra $A$, {\em i.e.}, a ring with a
ring homomorphism
$ \eta_\ahha : k \rightarrow Z(A)$ to
its centre. Denote by
$A^\mathrm{op}$ the
opposite algebra and by
$A^\mathrm{e} := A \otimes A^\mathrm{op}$
the enveloping algebra
of $A$, and by
$A-\Mod$ the category of
left $A$-modules.

The notions of  {\sl  $ A $--ring\/}  and  {\sl  $ A $--coring}  are direct generalizations of the notions of algebra and coalgebra over a commutative ring.

\begin{definition}  \label{def_A-coring}
 An  $ A $--coring  is a triple  $ (C , \Delta , \epsilon) $  where  $ C $  is an  $ A^e $--module  (with left action  $ L_A $  and right action  $ R_A $),  $ \; \Delta : C \longrightarrow C \otimes_A C \; $  and $ \; \epsilon : C \longrightarrow A \; $  are  $ A^e $--module  morphisms such that
  $$  (\Delta \otimes \text{\sl id}_C) \circ \Delta \; = \; (\text{\sl id}_C \otimes \Delta ) \circ \Delta  \quad  ,  \qquad  L_A \circ (\epsilon \otimes \text{\sl id}_C) \circ \Delta \; = \; \text{\sl id}_C \; = \; R_A \circ (\text{\sl id}_C \otimes \epsilon) \circ \Delta  $$
   \indent   As usual, we adopt Sweedler's  $ \Sigma $--notation  $ \, \Delta(c) = c_{(1)} \otimes c_{(2)} \, $  or  $ \, \Delta(c) = c^{(1)} \otimes c^{(2)} \, $  for  $ \, c \in C \, $.
   \end{definition}

 The notion of $A$-ring is dual to that of $A$-coring.  It is well known (see  \cite{Boe})  that  $ A $--rings  $ H $  correspond bijectively to  $ k $--algebra  homomorphisms  $ \, \iota : A \longrightarrow H \; $. A $A$-ring $H$ is endowed with an $A^e $-module structure :
 $$\forall h\in H, \quad a,b \in H, \quad a\cdot h \cdot b=\iota (a)h\iota(b).$$

\section{Preliminaries}
\label{tranquilli}

We recall the notions and results with respect to bialgebroids  that are needed to make  this article self content; see, {\em e.g.}, \cite{Kow} and references below  for an overview on this subject.

\subsection{Bialgebroids}
\label{h-Hopf_algbds}

For an  $ A^e $-ring $ U $  given by the $k$-algebra map
  $\eta: \Ae \to U $,  consider the restrictions
$ s := \eta( - \otimes 1_U) $  and  $t:= \eta(1_U \otimes -)$, called  {\it source\/}  and  {\it target\/}  map, respectively.
Thus an  $ A^e $-ring  $ U $ carries two  $ A $-module  structures from the left and two from the right, namely
  $$
a \lact u \ract b  :=  s(a)  t(b)  u,  \quad   \quad  a \blact u \bract b  :=  u  t(a)  s(b),   \eqno \forall \; a, b \in A  ,  u \in U.
$$
If we let  $  U_\ract {\otimes_{\scriptscriptstyle A}} {}_\lact U  $  be the corresponding tensor product of  $ U $  (as an $ A^e $-module)  with itself, we define the  {\it (left) Takeuchi-Sweedler product\/}  as
$$
U_\ract \! \times_\ahha \! {}_\lact U  \; := \;
     \big\{ {\textstyle \sum_i} u_i \otimes u'_i \in U_\ract \! \otimes_{\scriptscriptstyle A} \! {}_\lact U \mid {\textstyle \sum_i} (a \blact u_i) \otimes u'_i = {\textstyle \sum_i} u_i \otimes (u'_i \bract a), \ \forall a \in A \big\}.
$$
By construction,  $  U_\ract \! \times_{\scriptscriptstyle A} \! {}_\lact U  $  is an  $ \Ae $-submodule  of  $  U_\ract \! \otimes_{\scriptscriptstyle A} \! {}_\lact U  $;  it is also an  $ A^e $-ring via factorwise multiplication, with unit $  1_U \otimes 1_U  $  and  $ \eta_{{}_{U_\ract \times_{\scriptscriptstyle A} {}_\lact U}}(a \otimes \tilde{a}) := s(a) \otimes t(\tilde{a})$.

Symmetrically, one can consider the tensor product
$  U_\bract \otimes_\ahha \due U \blact {} $  and define the  {\em (right) Takeuchi-Sweedler product\/}  as
$U_\bract \times_\ahha \due U \blact {}  $,   which is an  $\Ae $-ring  inside
$  U_\bract \otimes_\ahha \due U \blact {} $.

\begin{definition}
 A {\em left  bialgebroid} $(U,A)$  is a  $ k $-module  $ U $  with the structure of an
$ \Ae $-ring  $(U, s^\ell, t^\ell)$  and an  $ A $-coring  $(U, \Delta_\ell, \epsilon )$  subject to the following compatibility relations:
\begin{enumerate}
\item
the  $ \Ae $-module  structure on the  $ A $-coring  $ U $  is that of
$ \due U \lact \ract  $;
\item
the coproduct $ \Delta_\ell $  is a unital  $ k $-algebra  morphism taking values in  $  U {}_\ract \! \times_{\scriptscriptstyle A} \! {}_\lact U  $;
\item
for all  $  a, b \in A  $,  $  u, u' \in U  $, one has:
\begin{equation}
\label{castelnuovo}
\epsilon (1_U)=1_A, \quad
\epsilon( a \lact u \ract b) =  a  \epsilon(u)  b, \quad \epsilon(uu')  =  \epsilon \big( u \bract \epsilon(u')\big) =  \epsilon \big(\epsilon(u') \blact u\big).
\end{equation}
\end{enumerate}
A  {\it morphism\/}  between left bialgebroids $(U, A)$ and $(U',A')$
is a pair $(F, f)$ of maps $F: U \to U'$, $f:A \to A'$ that commute with all structure maps in an obvious way.
\end{definition}

As for any ring, we can define the categories $U-\Mod$ and $\Mod -U$ of left and right modules over $U$. Note that 
$U-\Mod$ forms a monoidal category but $\Mod-U$ usually does not. However, in both cases there is a forgetful functor $U-Mod \to A^e-\Mod$, respectively  $\Mod-U \to A^e \Mod$ given by the formulas  : for $m \in M, \ n \in N, \ a,b \in A$
$$
a \lact m \ract b := s^\ell(a)t^\ell(b)m, \qquad a \blact m \bract b := ns^\ell(b)t^\ell(a)  
$$
For example, the base algebra $A$ itself is a left $U$-module via the left action
\begin{equation}\label{action of U on A}
u(a) := \epsilon( u \bract a) = \epsilon( a \blact u ), \quad  \forall u \in U, \quad   \forall a \in A  ,
\end{equation}
 but in most cases there is no right $U$-action on $A$.
 
 \begin{example}\label{the left bialgebroid V(L)}
 Let $L$ be a Lie Rinehart algebra (\cite{Rin}) over a commutative $k$-algebra $A$ with anchor $\rho : L \to Der_k (A)$. Its enveloping algebra $V(L)$ is endowed with a standard left bialgebroid described by : For all $a\in A$, $D \in L$, $u \in V(L)$,
 \begin{enumerate}
\item $s^\ell$ and $t^\ell$ are equal to the natural injection $\iota : A \to V(L)$
\item $\Delta_\ell: V(L) \to V(L)\otimes _A V(L), \quad \Delta (a)=a \otimes_A 1, \quad \Delta (D)=D \otimes_A 1+ 1\otimes_A D$
\item $\epsilon (u)=\rho (u)(1)$
\end{enumerate}
In this example, the left action of $V(L)$ on $A$ coincide with the anchor extended to $V(L).$. 
\end{example}

\subsection{Left and right Hopf left bialgebroids}
\label{goeseveron}
\label{half-Hopf_algebroids}

 For any  left  bialgebroid  $ U  $,  define the  {\em Hopf-Galois maps}
\begin{equation*}
\begin{array}{rclrcl}
\ga_\ell : \due U \blact {} \otimes_{\Aopp} U_\ract &\to& U_\ract  \otimes_\ahha  \due U \lact,
& u \otimes_\Aopp v  &\mapsto&  u_{(1)} \otimes_\ahha u_{(2)}  v, \\
\ga_r : U_{\!\bract}  \otimes^\ahha \! \due U \lact {}  &\to& U_{\!\ract}  \otimes_\ahha  \due U \lact,
&  u \otimes^\ahha v  &\mapsto&  u_{(1)}  v \otimes_\ahha u_{(2)}.
\end{array}
\end{equation*}

%
%

With the help of these maps, we make the following definition due to Schauenburg \cite{Sch}:

\begin{definition}
\label{def Half Hopf bialgebroids}
 A left bialgebroid $U$ is called a  {\em left Hopf left bialgebroid} or 
{\em $\times_A$ Hopf algebra}
if   $ \alpha_\ell  $ is a bijection.  Likewise, it is called a {\em right Hopf left bialgebroid} if $\ga_r$ is 
a bijection.
In either case, we adopt for all $u \in U$ the following (Sweedler-like)  notation
\begin{equation}
\label{latoconvalida}
u_+ \otimes_\Aopp u_-  :=  \alpha_\ell^{-1}(u \otimes_\ahha 1),  \quad
   u_{[+]} \otimes^\ahha u_{[-]}  :=  \alpha_r^{-1}(1 \otimes_\ahha u),
\end{equation}
and call both maps  $  u  \mapsto  u_+ \otimes_\Aopp u_-  $  and  $  u  \mapsto  u_{[+]} \otimes^\ahha u_{[-]}  $ {\em translation maps}.\\

\end{definition}


\begin{remarks}
\label{left/right-Hopf-left-bialg_cocomm}
\noindent 
Let $(U, A, s^\ell,t^\ell, \Delta , \epsilon)$ be a left bialgebroid. 
\begin{enumerate}

\item
In case $A=k$ is central in $U$, one can show that $\ga_\ell$ is invertible if and only if $U$ is a Hopf algebra, and the translation map reads
 $  u_+ \otimes u_-  :=  u_{(1)} \otimes S(u_{(2)})  $, where $S$ is the antipode of $U$.
On the other hand, $U$ is a Hopf algebra with invertible antipode if and only if both $\ga_\ell$ and $\ga_r$ are invertible,
and then $  u_{[+]} \otimes u_{[-]} := u_{(2)} \otimes S^{-1}(u_{(1)})  $.
%
\item
The underlying left bialgebroid in a {\em full} Hopf algebroid with bijective antipode is both a left and right Hopf left bialgebroid (but not necessarily vice versa); see \cite{BoeSzl} [Prop. 4.2] for the details of this construction recalled further.
\end{enumerate}

\end{remarks}

\begin{example} If $L$ is a Lie Rinehart algebra over a commutative $k$-algebra $A$ with anchor $\rho$, then its enveloping algebra $V(L)$, endowed with its standard bialgebroid structure, 
is a left  Hopf left bialgebroid. The translation map is  described as follows (in this case, $A=A^{op}$ and $s^\ell=t^\ell$) : If $a \in A$ and $D \in L$, 
$$a_+\otimes_{A^{op}} a_-=a \otimes_\Aopp 1, \quad D_+\otimes _\Aopp D_-=D\otimes_\Aopp  1- 1 \otimes_\Aopp D.  $$
 It is also a right Hopf algebroid as it is cocommutative. 
\end{example}

 The following proposition collects some properties of the translation maps  \cite{Sch}:

\begin{proposition}
Let $U$ be a left bialgebroid.
\begin{enumerate}

\item If  $  U $  is a left Hopf left bialgebroid, the following relations hold:
\begin{eqnarray}
\label{sch1}
u_+ \otimes_\Aopp  u_- & \in
& U \times_\Aopp U,  \\
\label{sch2}
u_{+(1)} \otimes_\ahha u_{+(2)} u_- &=& u \otimes_\ahha 1 \quad \in U_{\!\ract} \! \otimes_\ahha \! {}_\lact U,  \\
\label{sch3}
u_{(1)+} \otimes_\Aopp u_{(1)-} u_{(2)}  &=& u \otimes_\Aopp  1 \quad \in  {}_\blact U \! \otimes_\Aopp \! U_\ract,  \\
\label{sch4}
u_{+(1)} \otimes_\ahha u_{+(2)} \otimes_\Aopp  u_{-} &=& u_{(1)} \otimes_\ahha u_{(2)+} \otimes_\Aopp u_{(2)-},  \\
\label{sch5}
u_+ \otimes_\Aopp  u_{-(1)} \otimes_\ahha u_{-(2)} &=&
u_{++} \otimes_\Aopp u_- \otimes_\ahha u_{+-},  \\
\label{sch6}
(uv)_+ \otimes_\Aopp  (uv)_- &=& u_+v_+ \otimes_\Aopp v_-u_-,
\\
\label{sch7}
u_+u_- &=& s^\ell (\varepsilon (u)),  \\
\label{Sch8}
\varepsilon(u_-) \blact u_+  &=& u,  \\
\label{sch9}
(s^\ell (a) t^\ell (b))_+ \otimes_\Aopp  (s^\ell (a) t^\ell (b) )_-
&=& s^\ell (a) \otimes_\Aopp s^\ell (b),
\end{eqnarray}
where, in  \rmref{sch1},   we mean the Takeuchi-Sweedler product
\begin{equation*}
\label{petrarca}
   U \! \times_\Aopp \! U   :=
   \big\{ {\textstyle \sum_i} u_i \otimes v_i \in {}_\blact U  \otimes_\Aopp  U_{\!\ract} \mid {\textstyle \sum_i} u_i \ract a \otimes v_i = {\textstyle \sum_i} u_i \otimes a \blact v_i, \ \forall a \in A \big\}.
\end{equation*}
\item
There are similar relations for $u_{[+]}\otimes_A u_{[-]}$ if $  U $  is a right Hopf left bialgebroid (see \cite{CheGavKow} for an exhaustive list).
\end{enumerate}
\end{proposition}

The existence of a translation map if $U$ is  a left or right Hopf left bialgebroid makes it possible to endow $\Hom$-spaces and tensor products of $U$-modules with further natural $U$-module structures. These structures were systematically studied  in \cite{CheGavKow}( proposition 3.1.1). 

\begin{proposition}
\label{structures}
Let $(U, A)$ be a left bialgebroid, $M, M' \in U-\Mod$ and $N, N' \in \Mod-U$ be left resp.\ right $U$-modules, denoting the respective actions by juxtaposition.
\begin{enumerate}

\item
Let $(U,A)$ be additionally a left Hopf left bialgebroid.
\begin{enumerate}

\item
The $\Ae$-module $\Hom_{\Aopp}(M,M') $ carries a left  $ U $-module structure given by
\begin{equation}
\label{gianduiotto1}
(u\cdot f)(m)  :=  u_+ \big( f(u_-m) \big).
\end{equation}
\item
The $\Ae$-module $\Hom_\ahha(N,N')  $ carries a left $U$-module structure via
\begin{equation}
\label{lingotto1}
(u \cdot f)(n)  :=  \big( f(nu_+) \big)u_-.
\end{equation}
\item
The $\Ae$-module $\due N \blact {} \otimes_\Aopp M_\ract$
carries a right  $ U $-module  structure via
\begin{equation}
\label{superga1}
(n \otimes_\Aopp m)\cdot  u  :=  nu_+ \otimes_\Aopp u_-m.
\end{equation}
\end{enumerate}
\item
Let  $(U,A)$  be a right Hopf left bialgebroid instead.
\begin{enumerate}
\item
The $\Ae$-module  $  \Hom_\ahha(M,M')  $  carries a left  $ U $-module structure  given by
\begin{equation}
\label{gianduiotto2}
(u\cdot f)(m)  :=  u_{[+]}\big( f(u_{[-]}m) \big).
\end{equation}
\item
The $\Ae$-module
$\Hom_{\Aopp}(N,N')  $  carries a left  $ U $-module structure  given by
\begin{equation}
\label{lingotto2}
(u \cdot  f)(n)  :=  \big( f(nu_{[+]}) \big)u_{[-]}.
\end{equation}
\item
The $\Ae$-module
$ N_{\bract} \otimes^\ahha \due M \lact {}  $
carries a right  $ U $-module structure  given by
\begin{equation}
\label{superga2}
(n \otimes^\ahha m) \cdot u  :=  nu_{[+]} \otimes^\ahha u_{[-]}m.
\end{equation}
\end{enumerate}
\end{enumerate}
\end{proposition}

\begin{corollary} (\cite{CheGavKow}) Let $U$ be left and right left bialgebroid. 
The evaluation map
\begin{equation}
\label{lacartachenontagliaglialberi}
 P \bract {} \otimes_A \Hom_{A}(_\blact P  , _\blact N)_{\ract} \to N, \quad p  \otimes_A \phi  \mapsto  \phi(p)
\end{equation}
is a morphism of right $U$-modules for any $ N \in \Mod - U$.
\end{corollary}

\section{Poincaré duality}

We start by recalling the definition of a dualizing module introduced in \cite{CheGavKow}.

\begin{definition}
\label{dualizing module} (\cite{CheGavKow})
 Let $U$ be a left and right Hopf left bialgebroid over $A$.  
A right  $ U$-module  $ P $ satisfying the following conditions is called a dualizing module for $U$. 
\begin{enumerate}
\item The $A$-module $_\blact P$  is finitely generated projective;
\item
The left $U$-module morphism
 $$
A \to  \Hom_{A}(_\blact P,_\blact P) ,  \quad  a  \mapsto  \{ p \mapsto p \bract a \}
$$
 is an isomorphism;
\item
the evaluation map
\begin{equation}
\label{lacartachenontagliaglialberi}
 P \bract {} {\otimes_A}_\lact \Hom_{A}( _\blact  P, _\blact N) \to N, \quad p  \otimes_\Aopp  \phi  \mapsto  \phi(p)
\end{equation}
is an isomorphism for any $ N \in \Mod - U$.
\end{enumerate}
.
\end{definition}


We can now state twisted Poincaré duality:

\begin{theorem} \label{Poincare duality}
Let $U$ be a left and right Hopf left bialgebroid over $A$. Assume that 
\begin{enumerate}
\item $\Ext^{i}_U(A,U)=\{0\}\; if \; i\neq d$ and set $\Lambda =\Ext^{d}_U(A,U)$.   
\item The left $U$-module $A$ admits a finitely generated projective resolution of finite length
\item $A$ is endowed with a right $U$-module structure (denoted $A_R$)  for which it is a dualizing module 
\item Let ${\mathcal T}$ be the left $U$-module $Hom (_\blact {A_R }, _\blact \Lambda)$ (see proposition \ref{structures}). The $A$ module $_\lact {\mathcal T}$ and the 
$A^{op}$-module $ {\mathcal T}_\ract$ are projective.
\end{enumerate}
Then, for all $n \in \N$ there is an isomorphism
$$Ext^i_{U}(A,M)\simeq Tor_{U} \left ( A_R,  {\mathcal T}_\ract {\otimes_A} _\lact M \right ).$$ 
\end{theorem}

\begin{remark}\label{the case of V(L)} 
In the case where $U$ is the enveloping algebra $V(L)$ of a  finitely generated projective Lie Rinehart algebra $L$, the hypothesis are verified (see \cite{Che}).   More precisely, if $L$ is a projective  $A$-module with constant rank $n$, then   $\Ext^{i}_{V(L)}(A,V(L))=\{0\}$ if $i\neq n$.  
Moreover, $\Lambda^n (L^*)$ is endowed with a natural right $V(L)$-module structure determined as follows : 
$$\forall a \in A, \quad \forall D \in L, \quad \forall \omega \in \Lambda^n (L^*), \quad \omega \cdot a =a\omega, \quad \omega \cdot D=-L_D(\omega)$$
where $L_D$ is the Lie derivative with respect to $D$. 
Also the right $U$-modules  $\Ext^n _{V(L)}(A,V(L))$ and  $\Lambda^n (L^*)$ are isomorphic.
\end{remark}

{ \it Proof of theorem \ref{Poincare duality}:}

We will make use the following lemma : 

\begin{lemma} (\cite{KowKra}, lemma 16)
Let $U$ be a right Hopf left bialgebroid. Let  $N$ a right $U$-module,  $M$ and ${\mathcal T}$ two  left $U$-modules. One has an isomorphism of $k$-modules : 
$$(N_\bract{\otimes _A}_\lact {\mathcal T}){\otimes _U} M \simeq
N{\otimes _U}({\mathcal T}_{\ract}{\otimes _A}_\lact M )$$
\end{lemma}

Let $P^\bullet$ be a bounded finitely generated projective resolution of the left $U$-module $A$ and let $Q^\bullet$ be a  projective resolution of the left $U$-module $M$.  The following computation holds in $D^b(k-Mod)$, the bounded derived category of $k$-modules. 
$$\begin{array}{rcl}
RHom_{U}(A,M)&\simeq& Hom_{U}(P^\bullet,M)\\
                           & \simeq & Hom_{U}(P^\bullet ,U)\otimes_{U} M\\
                           &\simeq &\Lambda [-d] \otimes_{U} Q^\bullet\\
&\simeq & \left [ {A_R}_\bract {\otimes _A }_\lact {\mathcal T}  \right ]\otimes_{U} Q^\bullet\  [-d]\quad  (A_R \;{\rm  is\; a\; dualizing \;module})\\
&\simeq & {A_R}  {\otimes _U } \left ( {\mathcal T}_\ract {\otimes_A}_\lact Q^\bullet \right ) ({\rm previous \; lemma})\\
&\simeq &A_R \otimes^L_U  ({\mathcal T}_\ract {\otimes_A} _\lact  M)\\
\end{array}$$
The last isomorphism follows from the fact the $A$-module $_\lact {\mathcal T}$ is projective and from the following lemma
\begin{lemma} 
Denote by $^\ell U$ the left $U$-module structure on $U$ given by left multiplication. The map 
$$\begin{array}{rcl}
\alpha_r(M) : \; ^\ell U_\bract {\otimes _A}_\lact {\mathcal T} & \to & {\mathcal T}_{\ract}{\otimes _A}_\lact U\\
u \otimes t & \mapsto & u_{(1)} t \otimes u_{(2)}
\end{array}$$
is an isomorphism. One has $\alpha_r^{-1}(t\otimes u)=  u_{[+]}\otimes u_{[-]}t$.  Thus the $U$-module ${\mathcal T}_{\ract}{\otimes _A}_\lact U$ is projective if the $A$-module 
$_\lact {\mathcal T}$ is projective.$\Box$
\end{lemma}

\begin{remark} 
\begin{enumerate}
\item In the case where $U=A\otimes A^{op}$ (see examples \ref{examples Hopf algebroids}), 
$\Ext^i_U(A, M)$ is the Hochschild cohomology and we recover Van den Berg's Hochschild twisted Poincaré duality. Moreover, the beginning of the proof is similar to 
the proof of \cite{VDB} (theorem 1). 
\item The isomorphism $\Ext_{U}^n(A, ,M)\simeq Tor^{U}_{d-n} ( M_\ract {\otimes_A}_\blact \Lambda ,A )$ is proved in \cite{KowKra}. But, one can show that if the 
$A$-module $\Lambda_\bract$ is projective, $Tor^{U}_{d-n} ( M_\ract {\otimes_A}_\blact \Lambda ,A )\simeq Tor^{U}_{d-n} (   \Lambda ,M )$ and we get the beginning of our proof. 
\end{enumerate}
\end{remark}

In the case of full Hopf algebroids, there is a natural choice of right $U$-module structure on $A$.

\begin{free text}{\bf Reminder on full Hopf algebroids.}
\label{reminders_H-ads}
 Recall that a full Hopf algebroid structure (see, for example, \cite{Boe}) on a  $ k $-module  $ H $
consists of the following data:
\begin{enumerate}
\item
a left bialgebroid structure  $  H^\ell := ( H, A, s^\ell, t^\ell, \Delta_\ell  , \epsilon)  $
over a  $ k $-algebra  $ A  $;
\item
 a right bialgebroid structure  $  H^r := ( H, B, s^r, t^r, \Delta_r  , \partial)  $ over a  $ k $-algebra  $ B $;
\item
the assumption that the  $ k $-algebra  structures for  $ H $  in  {\em (i)\/}  and in  {\em (ii)\/} be the same;
\item
a $ k $-module  map  $  S : H \to H  $;
\item
some compatibility relations between the previously listed data for which we refer to {\em op.\ cit.}
\end{enumerate}
We shall denote by lower Sweedler indices the left coproduct  $ \Delta_\ell $  and by upper indices the right coproduct  $ \Delta_r  $, that is,
$  \Delta_\ell(h) =: h_{(1)} \otimes_\ahha h_{(2)}  $  and  $  \Delta_r(h) =: h^{(1)} \otimes_{B} h^{(2)}  $
 for any  $ h \in H  $.
A full Hopf algebroid (with bijective antipode) is both a left and right Hopf algebroid but not necessarily vice versa.
In this case, the translation maps in \rmref{latoconvalida} are given by
\begin{equation}
\label{laterza}
h_+ \otimes_\Aopp h_- = h^{(1)} \otimes_\Aopp S(h^{(2)}) \quad \mbox{and} \quad h_{[+]} \otimes_\Bop h_{[-]} = h^{(2)} \otimes_\Bop S^{-1}(h^{(1)}),
\end{equation}
formally similar as for Hopf algebras.
\end{free text}
The following lemma \cite{Boe, BoeSzl} will be needed to prove the main result in this subsection.

\begin{proposition}
\label{properties of S}
  Let  $ H =(H^\ell , H^r )$  be a (full) Hopf algebroid over $A$ with bijective antipode $S$.  Then
 \begin{enumerate}
\item

the maps $ \nu  :=  \partial s^\ell : A \to B^\op  $  and  $  \mu  :=  \epsilon s^r : B \to A^\op  $  are isomorphisms of  $ k $-algebras;
\item One has $\nu^{-1}=\epsilon t^r$ and $\mu^{-1}=\partial t^\ell$.
 \item
 the pair of maps $(S, \nu) : H^\ell \to {(H^r)}^\op_\coop  $  gives an isomorphism of left bialgebroids;
 \item
 the pair of maps
  $ (S, \mu) : H^r \to {(H^\ell)}^\op_\coop  $  gives an isomorphism of right bialgebroids.
 \end{enumerate}
 \end{proposition}

\begin{examples}\label{examples Hopf algebroids}
\begin{enumerate}
\item Let $A$ be a $k$-algebra, then $A^e=A \otimes A^{op}$ is a $A$-Hopf algebroid described by :
for all $a, b\in A$, 
$$\begin{array}{l}
s^\ell (a)=a \otimes 1, \quad t^\ell (b)=1\otimes b,\\
\Delta_\ell : A^e \to A^e \otimes_A A^e, \quad a\otimes b \mapsto (a\otimes_k 1)\otimes_A (1\otimes_k b)\\
\epsilon : A^e \to a , \quad a\otimes b \mapsto ab\\
s^r(a)=1 \otimes_k a, \quad t^r(b)=b\otimes_k 1\\
\Delta_r : A^e \to A^e \otimes_{A ^{op}}A^e, \quad a\otimes b \mapsto (1\otimes_k a)\otimes_A (b\otimes_k 1)\\
\partial : A^e \to a , \quad a\otimes b \mapsto ba\\
\end{array}$$

\item Let $A$  be a commutative  $k$-algebra and $L$ be a $k$-$A$-Lie Rinehart algebra.  Its enveloping algebra $V(L)$ is endowed with a standard left  bialgebroid structure
(see example \ref{the left bialgebroid V(L)}). Kowalzig (\cite{Kow}) has shown that  the left bialgebroid $V(L)$ can be endowed with a  Hopf algebroid structure if and only if there exists  a right $V(L)$-module structure on $A$. 
Then the right bialgebroid structure $V(L)_r$ is described as follows : for any $a\in A$, $D \in L$ and $u \in V(L)$: 
\begin{enumerate}
\item $\partial (u)=1 \cdot u$
\item $\Delta_r : V(L) \to V(L)_\bract {\otimes_A} _\blact V(L), \quad $
$\Delta_r (D)=D\otimes_A 1+1 \otimes_A D -\partial (X)\otimes_A 1$ and $\Delta_r(a)=a\otimes 1$
\item $S(a)=a $, $S(D)=-D+\partial (D)$. 
\end{enumerate}

It is in particular the case if $X$ is a ${\mathcal C}^\infty$ Poisson manifold, $A={\mathcal C}^\infty (X)$  and $L=\Omega^1(X)$ is the $A$-module of global differential 1- forms on $X$ . Huebschmann has shown (\cite{Hue})
that there is a $V(\Omega^1 (X))$-module stucture on $A$ determined by~: 
$$ a\cdot u=au  \quad {\rm and} \quad a\cdot udv=\{au,v\}.$$
Thus, $V(\Omega^1(X))$ is endowed with a (full) Hopf algebroid structure. 

\end{enumerate}
\end{examples}

 \begin{notation}
 Let $(H^\ell, H^r, S)$ be a full Hopf algebroid over $A$. 
 \begin{enumerate}
\item  Il $N$ is a right $H^\ell$ -module, we will denote by $_SN$ the left $H^\ell$-module defined by 
$$\forall h\in H,\quad  \forall n \in N, \quad  h\cdot_S n=n\cdot S(h).$$
\item Il $M$ is a left  $H^\ell$ -module, we will denote by $M_S$ the right  $H^\ell$-module defined by 
$$\forall h\in H, \quad \forall m \in M, \quad  m\cdot_S h=S(h)\cdot m.$$
\end{enumerate}
 \end{notation}
 
 \begin{remark}\label{structures of H-modules on A}
 If $H=(H^\ell, H^r, S)$ is a Hopf algebroid over a $k$-algebra $A$. We have the following module structures :

\begin{itemize}
\item a left $H^\ell$ -module structure given by $h\cdot_\ell a=\epsilon (hs^\ell(a))=\epsilon (ht^\ell(a))$.
\item a right $H^r$-module structure given by $\alpha\cdot_r h=\partial  (s^r(\alpha )h)=\partial  (t^r(a)h)$. 
\end{itemize}

Thanks to the proposition \ref{properties of S}, these two structures are linked by the relation 
$$S(h)\cdot_\ell \mu(\alpha ) = \mu [ \alpha  \cdot_r h].$$
\end{remark}

\begin{theorem} Let $(H^\ell, H^r)$ be a full Hopf algebroid over $A$ with bijective antipode $S$.  Consider $A$ with its left $H$-module structure (as in remark \ref{structures of H-modules on A}). 
We keep the notation of proposition \ref{properties of S}, in particular $\mu=\epsilon s^r$ and $\nu=\partial s^\ell$. 
 \begin{enumerate}
  \item If  $a \in A$, then $1 \cdot_S t^ \ell(a)=a$. Thus the $A$-module $_\blact{A_S}$ is free with basis 1. 
  \item If  $a \in A$, then $\alpha \cdot_S s^\ell(a)=\mu \nu (a) \alpha$. Thus  the $A^{op}$-module ${A_S}_\bract$ is free with basis 1. 
\item  If $N$ is a right $H^\ell$-module, the left $H^\ell$-module $Hom_A(_\blact {A_S} , _\blact N )$ is isomorphic to $_SN$. 
\item The right $H^\ell$-module $A_S$ is a dualizing module for $H$. 
\end{enumerate}
\end{theorem}

{\it Proof :}
\begin{enumerate}

\item  Using proposition \ref{properties of S}, we have

$$1 \cdot_S t^\ell (a)
= S(t^\ell (a))[1 ]
\underset{prop. \ref{properties of S}}{=} t^r\nu (a )[1]
=\epsilon  \left [ t^r\nu (a)\right ]=a
$$

\item Similarly, on has : $1\cdot_S s^\ell(a)=S(s^\ell (a))(1)=\epsilon s^r \nu (a) =\mu \nu (a).$


\item  The map 
$$\begin{array}{rcl}
Hom_A({A_S}_\bract , N_\bract )& \to & _SN\\
\lambda & \mapsto & \lambda (1)
\end{array}$$
is an isomorphism of left $H^\ell $-module as shows the following computation. Let $\alpha \in A_S$, $h \in H^\ell$ and 
$\lambda \in Hom_A(_\blact{A_S} , _\blact N)$. Using assertion 1. and theorem \ref{properties of S}, we have : 
$$\begin{array}{rcl}
(u\cdot \lambda )(1)&=&
\lambda(1 \cdot_S  h^{(1)})S(h^{(2)})\\
&=&\lambda \left [S(h^{(1)})(1)\right ]S(h^{(2)})\\
&=&\lambda \left [\epsilon S(h^{(1)})\right ]S(h^{(2)})\\
&=&\lambda \left [1\cdot_S t^\ell \epsilon S(h^{(1)})\right ]S(h^{(2)})\\
&=&\lambda (1) t^\ell \epsilon [S(h^{(1)})]S(h^{(2)})\\
&=&\lambda (1) t^\ell \epsilon [S(h)_{(2)})]S(h)_{(1)}\\
&=& \lambda (1) S(u)
\end{array}$$

\item  Let $N$ be a right $H^\ell$-module and let $n \in N$. Denote by $\lambda_n$ the element of $ \Hom_{A}( _\blact  {A_S}, _\blact N)$ determined by 
$\lambda_n (1)=n$.  By assertion 1 and  2,  
the map  $P \bract {} {\otimes_A}_\lact \Hom_{A}( _\blact  P, _\blact N) \to N, \quad p  \otimes_\Aopp  \phi  \mapsto  \phi(p)$ is an isomorphism with inverse 
$n \mapsto 1\otimes \lambda_n$. 

We need now to check that the map  $ A  \to  \Hom_{A}(_\blact {A_S},_\blact  {A_S}) ,  \quad  a  \mapsto  \{ p \mapsto p \bract a \}$ is an isomorphism.
By assertion 3, this boils down to showing that 
$ A  \to {}_S(A_ S), \quad  a  \mapsto  1 \bract a $ is an isomorphism .
But, this is true as $1 \bract a= \mu \nu (a)$. Indeed, 
$$1\bract a =S^2(s^\ell (a))(1)= \epsilon S^2\left [ s^\ell(a)\right ]=\mu \partial \left [ S(s^\ell (a)) \right ]=\mu \nu \epsilon (s^\ell (a))=\mu \nu (a). \Box$$

\end{enumerate}

We can now state  twisted Poincaré duality for Hopf algebroids.

\begin{theorem} Let $(A, H^\ell, H^r)$ be a Hopf algebroid over $A$ with bijective antipode $S$.  As in proposition \ref{properties of S}, we set $\mu=\epsilon s^r$ and $\nu= \partial s^\ell$. Assume that 
\begin{itemize}
\item $\Ext^{i}_{H^\ell}(A,H^\ell)=\{0\}\; if i\neq d$ and set $\Lambda =\Ext^{d}_{H^\ell}(A,H^\ell)$.  
\item $_\blact \Ext^{d}_{H^\ell}(A,H^\ell)$ is a   projective $A$-module and 
 $ \Ext^{d}_{H^\ell}(A,H^\ell)_\bract $ is a   projective $A^{op}$-module
\item The left $H^\ell$-module $A$ admits a finitely generated projective resolution of finite length  
\end{itemize}
Then there is an isomorphism
$$\Ext^i_{H^\ell}(A,M)\simeq \Tor_{d-i}^{H^\ell} \left ( A_S,  {_S\Lambda}_\ract {\otimes_A} _\lact M \right )  $$ 
\end{theorem}

As an application, we find a Poincaré duality for smooth Poisson algebras. Assume that $X$ is a ${\mathcal C}^\infty (X)$ Poisson manifold, $L=\Omega^1(X)$, 
Huebschmann has shown that the $\Ext_{V(\Omega^1(X))}^i (A, M)$'s  coincide with the Poisson cohomology with coefficients  in $M$, $H^i _{Pois}(A,M)$. 
Also, the $\Tor_i^A(A_S, M)$'s  coincide with the Poisson homology with coefficients in $M$, $H_i^{Pois}(A,M)$. 

\begin{corollary}
Let $X$ be a ${\mathcal C}^\infty$ $n$-dimensional Poisson manifold, $A={\mathcal C}^\infty (X)$ and $M$ a left $V(\Omega^1(M))$-module. 
Let $S$ be the antipode of the (full) Hopf algebroid $V(\Omega^1(X))$ (see examples \ref{examples Hopf algebroids}).
Then 
${\mathcal T}$ is isomorphic to  $_S\left ( \Lambda^n \Omega^1(X)^*\right  )=_S\Lambda^n Der(A)$ where $df$ acts (on the right) as the opposite of the Lie derivative of the Hamiltonian vector field $X_f$. There is an isomorphism 
$$H^i_{Pois}(A, M)
\simeq H_{d-i}^{Pois}(A, \; _S{\Lambda^{n}  Der (A)}   \otimes_A M ).$$
\end{corollary}

\begin{remark}
This formula is proved  in  \cite{CLYZ} for oriented Poisson manifolds and $M=A$ (see also \cite{LauRic} for polynomial algebras with quadratic Poisson structures, \cite{Zhu} for linear Poisson structures, \cite{LuoWangWu} for general polynomial Poisson algebras)
\end{remark}

Sophie Chemla, Sorbonne Université, sophie.chemla@sorbonne-universite.fr

\end{document}